\documentclass[preprints,article,accept,moreauthors,pdftex]{Definitions/mdpi} 

\usepackage{bm}
\usepackage{bbm}
\usepackage{listings}
\usepackage{subfig }
\usepackage{amsmath}
\usepackage{amssymb}
\usepackage{amsthm}
\usepackage{booktabs}
\usepackage{hyperref}
\usepackage{erewhon} % serif
\usepackage[sf,scale=1]{libertine}
\usepackage{fourier}
\usepackage[T1]{fontenc}
\usepackage{algorithm}
\usepackage{algorithmic}

\usepackage{tikz}
\usepackage{placeins}
\usetikzlibrary{decorations.pathreplacing,angles,quotes}
\usetikzlibrary{arrows}
\usetikzlibrary{arrows,shapes, positioning,fit,calc, decorations.pathmorphing,matrix}
\tikzstyle{inputNode}=[draw=none,rectangle, rounded corners, minimum size=10pt,inner sep=2pt,fill=green!50]
\tikzstyle{stateTransition}=[->, decorate, thick, shorten <= 2pt, shorten >= 2pt]
\tikzstyle{hiddenA}=[draw,rectangle, text width=360pt,	minimum width=360pt, minimum height=12pt,inner sep=2pt,fill=red!10]
\tikzstyle{hiddenB}=[draw,rectangle, text width=360pt,	minimum width=360pt, minimum height=12pt,inner sep=2pt,fill=yellow!20]
\tikzstyle{hiddenC}=[draw,rectangle, text width=360pt,	minimum width=360pt, minimum height=12pt,inner sep=2pt,fill=magenta!15]
\tikzstyle{output}=[draw=none,fill=blue!30,rectangle, rounded corners,minimum size=10pt,inner sep=1pt]

\hypersetup{colorlinks,citecolor=teal}

\firstpage{1} 
\pubvolume{xx}
\issuenum{1}
\articlenumber{1}
\pubyear{2022}
\copyrightyear{2022}
\history{}

\newcommand{\arcsinh}{\mathop{\mathrm{asinh}}} 

\newcommand*\diff{\mathop{}\!\kern0pt\mathrm{d}}

\Title{Inserting or Stretching Points in Finite Difference Discretizations}
\Author{Jherek Healy}
\AuthorNames{Jherek Healy}
\address{}
\corres{Correspondence: jherekhealy@protonmail.com}
\abstract{Partial differential equations sometimes have critical points where the solution or some of its derivatives are discontinuous. The simplest example is a discontinuity in the initial condition. It is well known that those decrease the accuracy of finite difference methods. A common remedy is to stretch the grid, such that many more grid points are present near the critical points, and fewer where the solution is deemed smooth. An alternative solution is to insert points such that the discontinuities fall in the middle of two grid points. This paper compares the accuracy of both approaches in the context of the pricing of financial derivative contracts in the Black-Scholes model and proposes a new fast and simple stretching function.}
\keyword{finite difference method, grid stretching, Black-Scholes}

\begin{document}

\section{Introduction}

Partial differential equations (PDEs) sometimes have critical points where the solution or some of its derivatives are discontinuous. The simplest example is a discontinuity in the initial condition. This situation arises in the pricing of nearly all financial derivative contracts. The vanilla European option of given maturity and strike price, the simplest non-linear contract, has indeed a discontinuous first derivative at the strike price.

It is well known that such critical points decrease the accuracy of finite difference methods. A common remedy, detailed in \citep[p. 167]{TaRa00}, is to stretch the grid such that many more grid points are present near the critical points, and fewer where the solution is deemed smooth. The stretching transformation for a single point reads
\begin{equation}
	S(u) = B + \alpha \sinh\left(c_2 u + c_1 (1-u)\right)\,,
\end{equation}
where $c_1 = \arcsinh\frac{S_{\min} - B}{\alpha}$, $c_2 = \arcsinh\frac{S_{\max} - B}{\alpha}$, and $\alpha$ controls the density of points near the critical point $B$. For $u \in [0,1]$, we have $S(u) \in [S_{\min},S_{\max}]$.

Independently of such a stretching, \citet{TaRa00, giles2005convergence} also show that the error in the solution is significantly decreased when the critical points are located in the middle of two grid points. There are several ways to place the critical points in such manner. A first approach is to move the grid. This is applicable only for a single critical point, and if the boundaries can be moved. A second approach is to simply insert a point in the grid, around the critical point such that the critical point is exactly in the middle of two grid points. A third approach is to use a smooth deformation, typically a monotonic cubic spline, to place the critical point approximately (but not exactly) in the middle of two grid points \citep[p. 171]{TaRa00}.

The advantage of the cubic spline smooth deformation is to preserve the second-order convergence. A robust implementation is however more involved than the insertion approach. The insertion approach, due to its lack of smoothness, will a priori not preserve the second-order convergence, but this does not mean that its accuracy is worse.

In this paper, we compare the accuracy of the two approaches, using concrete examples of options in the Black-Scholes model, on nearly uniform grids, as well as on stretched grids. We also propose a faster stretching transformation, similar to the sinh transformation and give a simple extension to multiple critical points.

\section{Cubic stretching}
\subsection{Single critical point}
According to \citet[p. 307]{Noye1983}, a stretching function should have the following properties:
\begin{enumerate}[label=(\roman*)]
	\item $\diff S/\diff u$ should be finite over the whole interval - if it becomes infinite at some point, then there is poor resolution near that point;
	\item $\diff S/\diff u$ must be smaller near at critical point than elsewhere in the interval, which ensures high resolution near the critical point, but $\diff S/\diff u$  should be non zero at the critical point.
\end{enumerate}

An intuitive candidate would be a function based on a probability density function. A mixture distribution makes it easy to ensure a higher density around the critical points. A numerical inversion of the mixture distribution, for example via a monotonic interpolation scheme, leads to the desired stretching function. Unfortunately, such a stretching will typically have very large derivatives near the boundaries (corresponding to the inverse of the cumulative density tails) and thus does not obey property (i).

For a single critical point, an interesting stretching function candidate is the cubic based on the Taylor series of the sinh function:
\begin{equation}
	S(u) = B + \alpha \left[\frac{1}{\chi}\left(c_2 u +c_1(1-u)\right)^3 + c_2 u +c_1(1-u)\right]\,,\label{eqn:sinh_cubic}
\end{equation}
where $c_1$ is the solution of the depressed cubic equation $\frac{1}{\chi}c_1^3 + c_1 + \frac{B-S_{\min}}{\alpha}=0$ and $c_2$ is the solution of $\frac{1}{\chi}c_2^3 + c_2 + \frac{B-S_{\max}}{\alpha}=0$. The value $\chi=6$ matches the sinh expansion, other positive values are also possible.
\begin{figure}[h]
	\begin{center}
\includegraphics[width=0.9\textwidth]{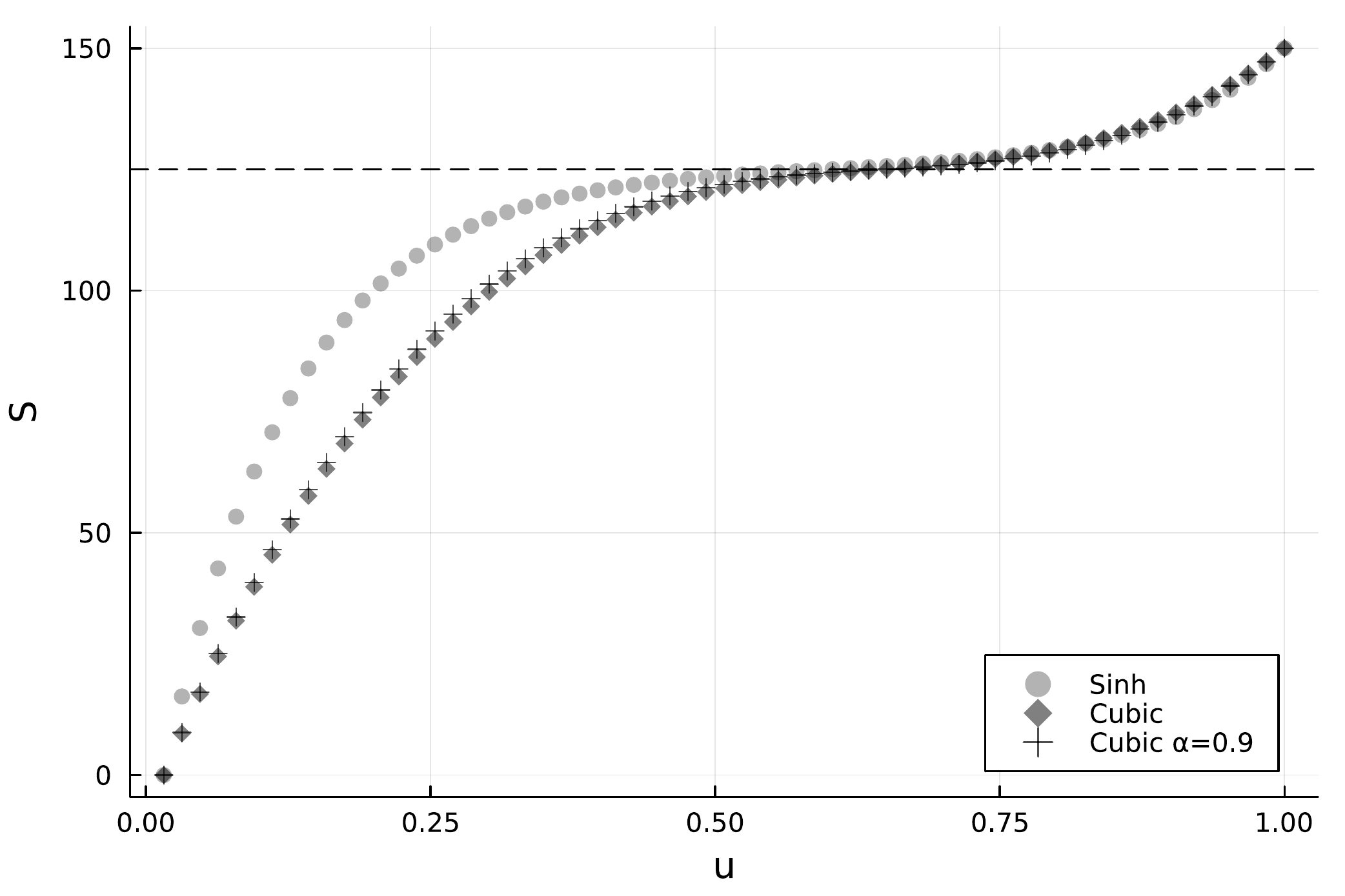}
	\end{center}
	\caption{Stretching around the point $B=125$ using 63 points in the interval $[0, 150]$ with $\alpha = 1.50$. \label{fig:sinh_cubic_points}}
\end{figure}

Figure \ref{fig:sinh_cubic_points} shows the cubic transformation to be close to the sinh transformation in practice. As expected, it is not exponential and thus closer to linear, far away from the critical point. For the same value of $\alpha$, the slope is slightly different at the critical point. The slope is matched using a lower $\alpha=0.9$ for the cubic stretching. One main advantage of the cubic stretching is performance, as the transformation doesn't involve any costly function at all. In practice, the cubic stretching is around five times faster.

\subsection{Many critical points}
\citet{TaRa00} propose to use the following jacobian for multiple critical points $(B_k)$:
\begin{equation}
	J(u,S) = \frac{\partial S}{\partial u} = A \left(\sum_{k} \frac{1}{\alpha_k^2+ (S(u)-B_k)^2}\right)^{-\frac{1}{2}}\label{eqn:hyperbolic_jacobian}\,,
\end{equation}
where $A$ is a normalizing constant used to ensure that $S(1) = S_{\max}$ with initial condition $S(0) = S_{\min}$. The Jacobian is nearly constant for $S \approx B_k$ which corresponds to a uniform discretization and is nearly linear for $S \gg B_k$ or $S \ll B_k$ which corresponds to a exponential grid.

Equation \ref{eqn:hyperbolic_jacobian} is an ordinary differential equation (ODE) whose initial condition consists in the function values at two end-points: it is a two-points boundary problem. A standard method to solve this kind of problem is  the shooting method: we are shooting a projectile from point $S(0) = S_{\min}$ so that it lands at point  $S(1) = S_{\max}$. Any solver can be used so solve for $A$. The ODE can be solved with the fourth-order Runge-Kutta method for a given guess $A$.

Similarly, the derivative of Equation \ref{eqn:sinh_cubic} provides a candidate stretching for multiple points:
\begin{equation*}
	\frac{\diff S}{\diff u} = \alpha A \prod_{i=1}^{n} (u-b_i)^2 + \alpha\,.
\end{equation*}
The solution $(A, b_1,...,b_n)$ such that $S(0) = S_{\min}$, $S(1) = S_{\max}$ and $S(b_i)=B_i$ involves a $n$-dimensional non-linear optimization and may not be practical for large $n$.

Solving such non-linear problems makes the overall technique much slower than the single critical point case, and more challenging to implement in a robust fashion. We thus describe below simpler, better performing and more robust techniques below.
\subsubsection{Direct piecewise-cubic representation}
Based on Equation \ref{eqn:sinh_cubic}, we consider a piecewise-cubic representation of class $\mathcal{C}^1$. Let $(B_i)_{i=1,...,m}$ be the ordered $m$ critical points in the interval $(S_{\min},S_{\max})$. Let $D_i = \frac{B_{i}+B_{i+1}}{2}$ be the corresponding mid-points for $i = 1,...,m-1$, and $D_0 = S_{\min}, D_{m}=S_{\max}$ for notation convenience. The piecewise cubic interpolant on the interval $[d_{i-1},d_i)$ reads
\begin{subequations}
	\begin{equation}
		p_i(u) = B_i +  \alpha_i \left[\frac{1}{\chi}\left(c_{2i} (u-d_{i-1}) +c_{2i-1}(d_i-u)\right)^3 + c_{2i} (u-d_{i-1}) +c_{2i-1}(d_i-u)\right]\,,\label{eqn:sinh_cubic_i} 
	\end{equation}
where $d_i$ is such that
   \begin{equation}
   	p_i(d_i) = D_i\,, \quad p_i(d_{i-1}) = D_{i-1}\,.
   	\end{equation}
   In particular we have $d_0 = 0$ and $d_m = 1$.
\end{subequations}
The variables $(c_i)_{i=1,...,m}$ are thus solutions of the following cubic equations
\begin{subequations}
	\begin{equation}
	\frac{1}{\chi}\left(c_{2i} (d_i-d_{i-1})\right)^3 + c_{2i} (d_i-d_{i-1}) +\frac{B_i-D_i}{\alpha_i} = 0\,,
\end{equation}
	\begin{equation}
	\frac{1}{\chi}\left(c_{2i-1} (d_i-d_{i-1})\right)^3 + c_{2i-1} (d_i-d_{i-1}) +\frac{B_i-D_{i-1}}{\alpha_i} = 0\,.
\end{equation}
This leads to the values $c'_{2i} = c_{2i}  (d_i-d_{i-1})$ and $c_{2i-1}'= c_{2i-1}(d_i-d_{i-1})$  for $i=1,...m$.
Furthermore, the continuity of the first derivative at $d_i$ imposes
	\begin{equation}
 \alpha_i \left[	\frac{3}{\chi}(c_{2i}-c_{2i-1}) \left(c_{2i} (d_i-d_{i-1})\right)^2 + c_{2i}-c_{2i-1}\right]= 	 \alpha_{i+1}\left[\frac{3}{\chi}(c_{2i+2}-c_{2i+1}) \left(c_{2i+1} (d_{i+1}-d_{i})\right)^2 +c_{2i+2}-c_{2i+1}\right]\,,
\end{equation}
or equivalently
	\begin{equation*}
 \alpha_i \left[	\frac{3}{\chi}(c_{2i}-c_{2i-1}) c_{2i}'^2 + c_{2i}-c_{2i-1}\right]= \alpha_{i+1}\left[	\frac{3}{\chi}(c_{2i+2}-c_{2i+1})c_{2i+1}'^2 +c_{2i+2}-c_{2i+1}\right]\,.
\end{equation*}
Multiplying by $(d_i-d_{i-1})(d_{i+1}-d_{i})$ leads to the tridiagonal system
	\begin{equation*}
	\alpha_i \left[\frac{3}{\chi}(c_{2i}'-c_{2i-1}') c_{2i}'^2 + c'_{2i}-c'_{2i-1}\right](d_{i+1}-d_{i})= 	\alpha_{i+1} \left[\frac{3}{\chi}(c_{2i+2}'-c_{2i+1}')c_{2i+1}'^2 +c_{2i+2}'-c_{2i+1}' \right](d_i-d_{i-1})\,,
\end{equation*}
which gives $d_i$ for $i=2,...,m-1$. From $d_i$ and $c_i'$, we trivially deduce the coefficients $c_i$.
\end{subequations}

The second derivative at each critical point is discontinuous, in fact it can be shown that, for a constant $\alpha_i=\alpha$, we have $p_{i}''(d_i)=-p_{i+1}''(d_i)$. Note that changing $\chi$ or moving $D_i$ will not help with the discontinuity. The discontinuity in the second derivative at the critical points is a fundamental aspect of our choice of piecewise-cubic representation, as the first derivative on each segment is a parabola with positive curvature. Intermediate knots would be required to derive a $\mathcal{C}^2$ interpolant.

The value of $\alpha$ may be made dependent on the critical point. For example we may choose a larger $\alpha$ for points corresponding to discontinuous second derivative in the solution, compared to points corresponding to a discontinuous first derivative in the solution.%TODO plot first derivative? or first difference?
\begin{figure}[h!]
	\begin{center}
				\subfloat[][Stretching function.]{\includegraphics[width=0.5\textwidth]{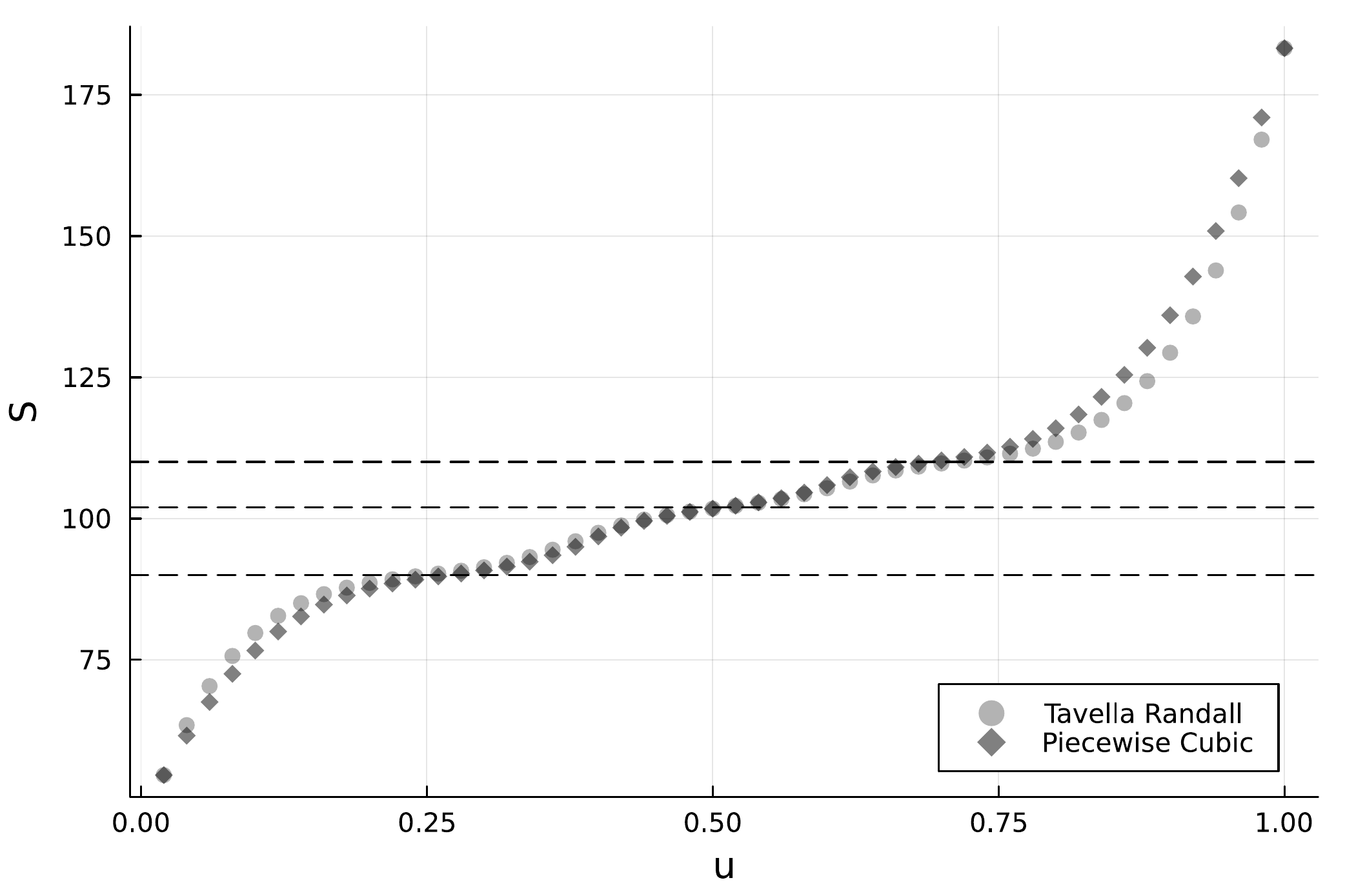}}
				\subfloat[][Derivative by forward difference.]{\includegraphics[width=0.5\textwidth]{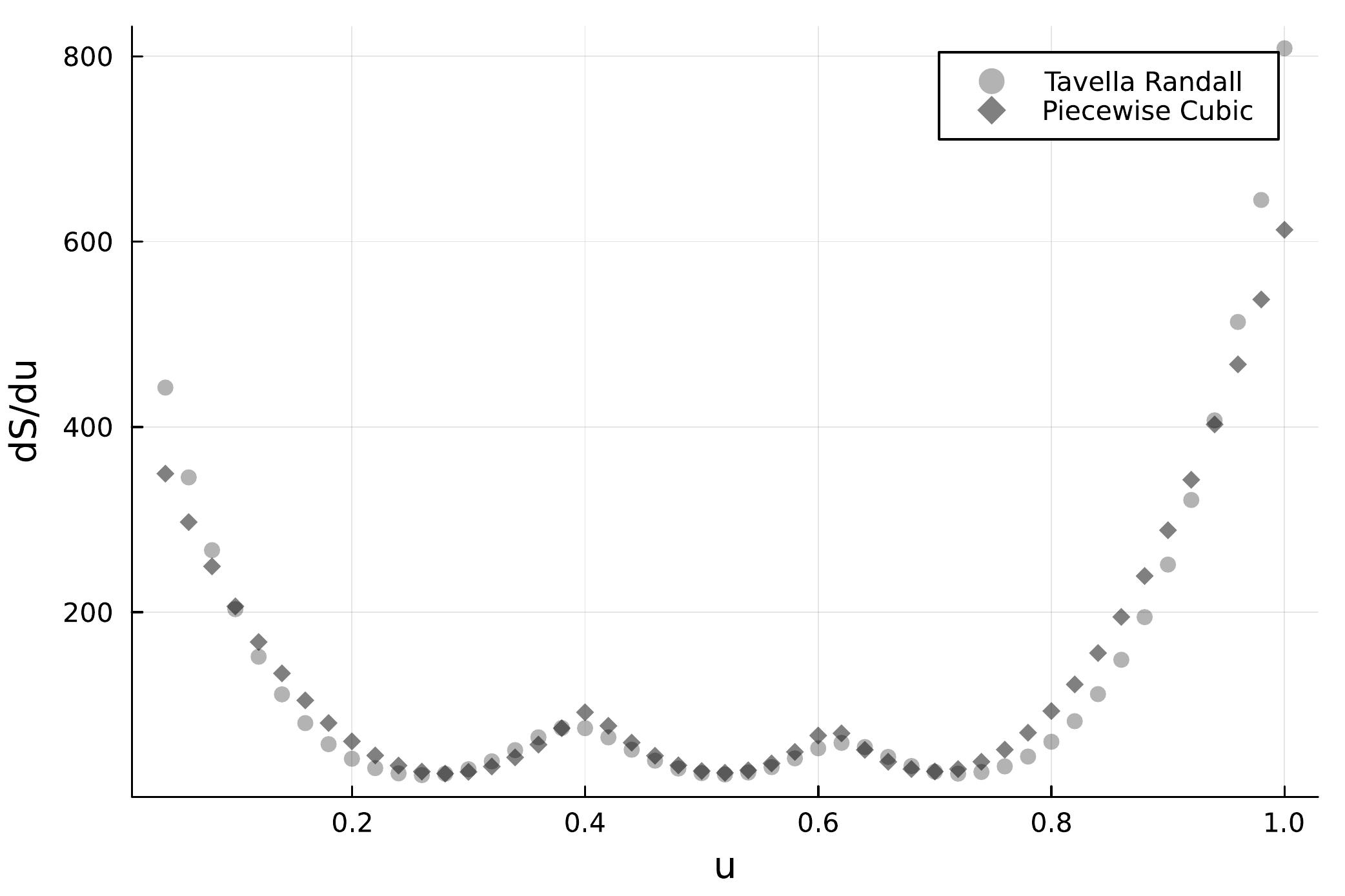}}
	\end{center}
	\caption{Stretching around the points $B_1=90, B_2=102, B_3=110$,  using 50 points in the interval $[54, 183]$ with $\alpha = 1.30$. \label{fig:sinh_cubic_points_multi}}
\end{figure}

\subsubsection{$\mathcal{C}^2$ piecewise representation}
It is possible to fill the piecewise cubic around the discontinuity in the second-derivative by a quintic piece  $\hat{p}_i(u)= \sum_{k=0}^5 a_k (u-d_i^L)^k$ such that
\begin{align*}
	\hat{p}_i(d_i^L) = p_i(d_i^L)\,,\quad 	\hat{p}_i'(d_i^L) &= p_i'(d_i^L)\,, \quad	\hat{p}_i''(d_i^L) = p_i''(d_i^L)\,,\\
	\hat{p}_i(d_i^R) = p_{i+1}(d_i^R)\,,\quad 	\hat{p}_{i}'(d_i^R) &= p_{i+1}'(d_i^R)\,,\quad 	\hat{p}_i''(d_i^R) = p_{i+1}''(d_i^R)\,.
\end{align*}
There is some flexibility towards the choice of $d_i^L$ and $d_i^R$ such that $d_i^L < d_i < d_i^R < d_{i+1}$.

A first candidate is to use only the $(d_i)$:
\begin{equation} d_i^L = d_i - \lambda \left(d_i-d_{i-1}\right)\,,\quad d_i^R = d_i + \lambda \left(d_{i+1}-d_{i}\right)\,,\label{eqn:quintic_knot_direct} \end{equation} with $\lambda \leq 1/2$. Typically, we pick $\lambda= 1/4$.

The inverse of the critical point is not necessarily in the middle of two $d_i$, and thus it may make more sense to use instead
\begin{equation}
	d_i^L = p_i^{-1}\left(B_{i-1} + (B_i-B_{i-1})\left(\frac{1}{2}-\lambda\right)\right)\,,\quad 	d_i^R = p_{i+1}^{-1}\left(B_{i} - (B_i-B_{i-1})\left(\frac{1}{2}-\lambda\right)\right)\,.\label{eqn:quintic_knot_inverse}
\end{equation}

The quintic coefficients $a_i$ thus obey
\begin{subequations}
\begin{equation}
	a_0 =  p_i(d_i^L)\,,\quad  a_1 = p_i'(d_i^L)\,,\quad a_2 = \frac{1}{2} p_i''(d_i^L)\,,
\end{equation}
\begin{equation}
	\begin{pmatrix}
		 \Delta_i^3 & \Delta_i^4 & \Delta_i^5\\
		 3\Delta_i^2 & 4\Delta_i^3 & 5\Delta_i^4\\
		 6\Delta_i & 12 \Delta_i^2 & 20 \Delta_i^3\\	 
	\end{pmatrix} \cdot 	\begin{pmatrix} a_3 \\ a_4 \\ a_5
\end{pmatrix}  =  	\begin{pmatrix} p_{i+1}(d_i^R) - (a_0 + a_1 \Delta_i + a_2 \Delta_i^2) \\
p_{i+1}'(d_i^R) - (a_1 + 2 a_2 \Delta_i)\\
p_{i+1}''(d_i^R) - 2 a_2
\end{pmatrix}
\end{equation}	
\end{subequations}
with $\Delta_i = (d_i^R - d_i^L)$. 

Furthermore, we want the quintic to be monotonic. This is achieved if the roots of $\hat{p}_i''$ in the interval $(d_i^L, d_i^R)$ are such that $\hat{p}_i' > 0$. A minor remark: we conjecture that a sufficient condition for monotonicity is $a_3 >0$. If the condition is not verified, we revert to the $\mathcal{C}^1$ representation. 

\begin{figure}[h!]
	\begin{center}
		\subfloat[][First derivative.]{\includegraphics[width=0.5\textwidth]{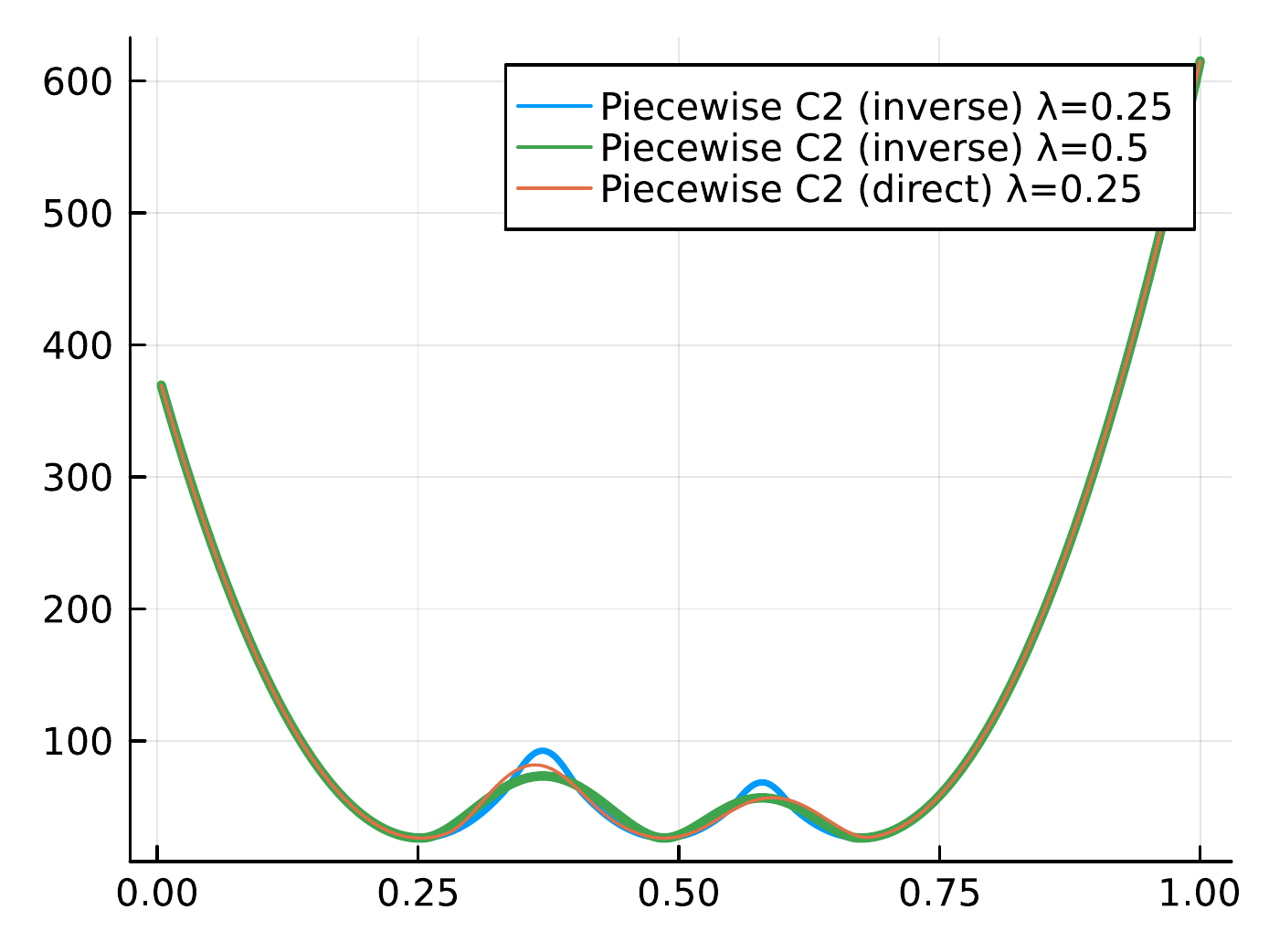}}
		\subfloat[][Second derivative.]{\includegraphics[width=0.5\textwidth]{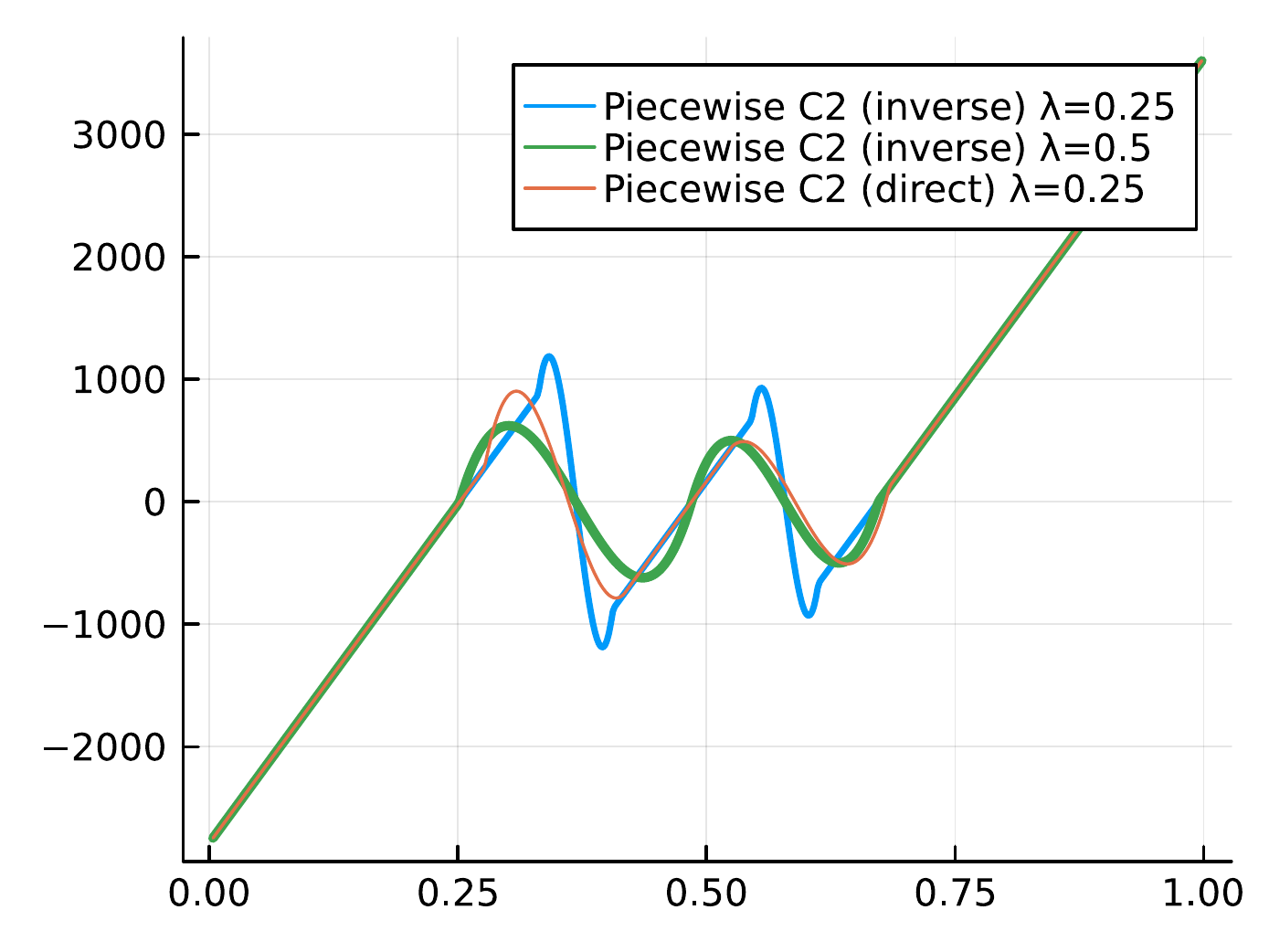}}
	\end{center}
	\caption{$\mathcal{C}^2$ Stretching around the points $B_1=90, B_2=102, B_3=110$ in the interval $[54, 183]$ with $\alpha = 1.30$. \label{fig:c2_cubic_points_d2l}}
\end{figure}

Figure \ref{fig:c2_cubic_points_d2l} shows little difference in the first derivative $\hat{p}_i'$ between the different choices. In particular, the simpler direct choice (Equation \ref{eqn:quintic_knot_direct}) is not necessarily worse. For the inverse approach  (Equation \ref{eqn:quintic_knot_inverse}), a large $\lambda=1/2$, which corresponds to starting the quintic at the critical points, leads to a smoother second derivative and looks surprisingly acceptable.

\begin{figure}[h!]
	\begin{center}
		\includegraphics[width=0.9\textwidth]{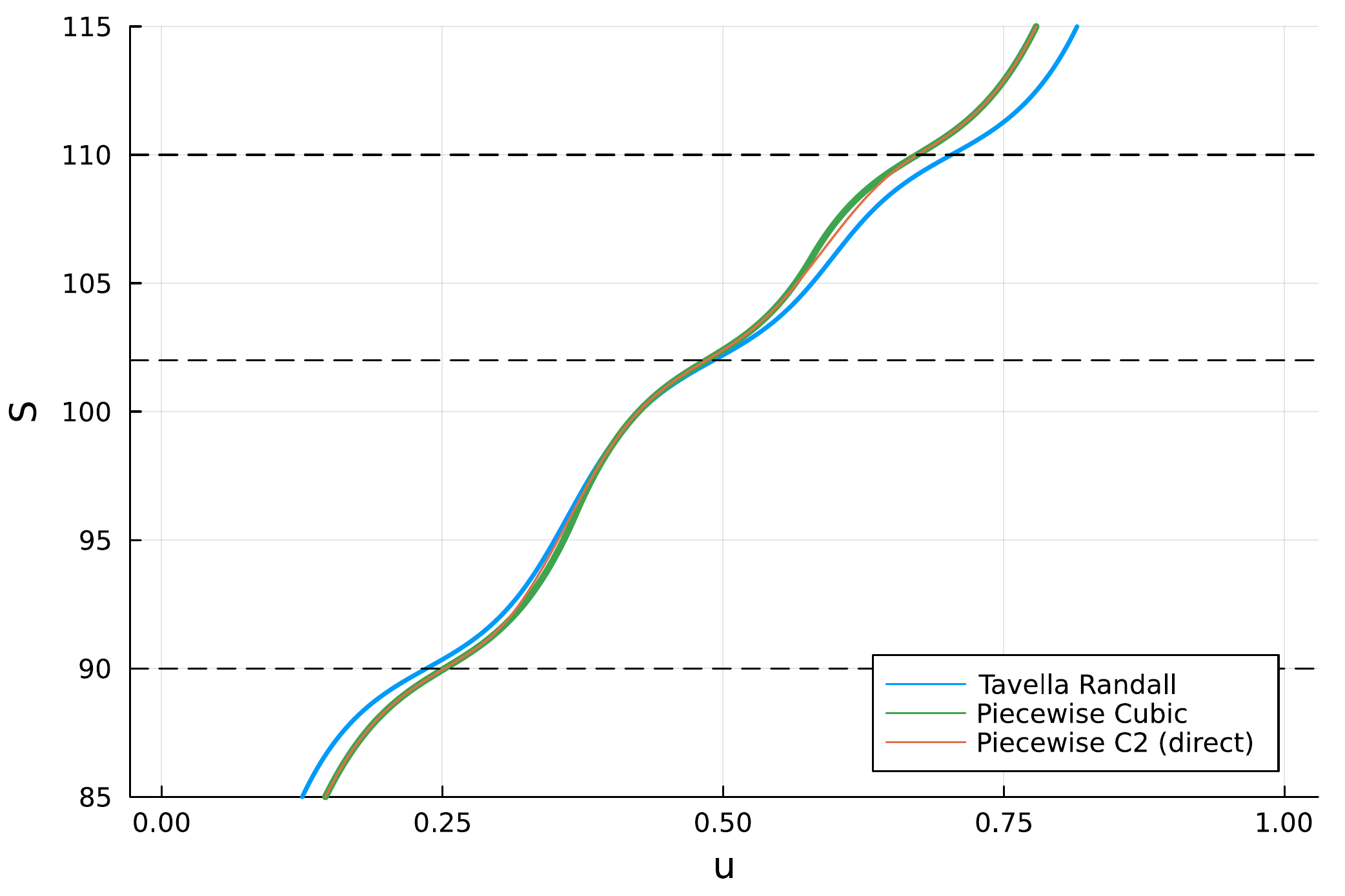}
	\end{center}
	\caption{Stretching around the points $B_1=90, B_2=102, B_3=110$ in the interval $[54, 183]$ with $\alpha = 1.30$. \label{fig:c2_cubic_points}}
\end{figure}

The stretched grid points are almost indistinguishable between the different piece-wise stretchings. Figure \ref{fig:c2_cubic_points} shows that the $\mathcal{C}^1$ cubic stretching leads to virtually the same points on this example.
%% PLOT lambda=0.25, lambda=0.5 for pinv based quintic second der, first der, func
%% PLOT lambda=0.25 for d based quintic => d looks smoother/closer to sinh. It seems ok to have lambda close to 0.5 for pinv and thus even if inverse is close to split, it's not such a big deal a priori.
%% The various choices do not seem to impact the stretched points significantly.
%% the cubic second derivative goes through 0=min, then 0=max then 0=min.
 % => if f'''(bi) > 0 => monotonic. We know c2-c1 > 0 / cubic has positive coeff because of 5a-5b. So true for the cubic but not obvious for the quintic.
%c.c.t.analytics.finitediff.lattice.transform.PointDuplicatePerformance.benchmarkInsertPointsArray   avgt   20    30.976 ±  15.255  us/op
%c.c.t.analytics.finitediff.lattice.transform.PointDuplicatePerformance.benchmarkInsertPointsNearly  avgt   20    17.662 ±   4.467  us/op
%c.c.t.analytics.finitediff.lattice.transform.PointDuplicatePerformance.benchmarkInsertPointsSinh    avgt   20   146.876 ±  21.131  us/op
%c.c.t.analytics.finitediff.lattice.transform.PointDuplicatePerformance.benchmarkInsertPointsTree    avgt   20   164.658 ±  33.997  us/op
%c.c.t.analytics.finitediff.lattice.transform.PointDuplicatePerformance.benchmarkNPointsPolynomial   avgt   20    20.169 ±   4.873  us/op
%c.c.t.analytics.finitediff.lattice.transform.PointDuplicatePerformance.benchmarkNPointsSinh         avgt   20  2228.282 ± 408.804  us/op
%c.c.t.analytics.finitediff.lattice.transform.PointDuplicatePerformance.benchmarkPointsPolynomial    avgt   20    15.292 ±   4.271  us/op

\section{Numerical Results}
We consider the same knock-out barrier option of maturity $T=1$ year, strike $K=100$ and barrier $B=125$,  with 250 discrete observations dates, starting at $t_1=1/250$ until $t_{250}=T=1$ under the Black-Scholes model with dividend yield $q=0.02$, interest rate $r=0.07$ and volatility $\sigma=20\%$, presented in \citep[Tables 6.1 and 6.2]{TaRa00}. The grid boundaries are at $S_{\min} ) = 0$ and $S_{\max} = 150$.

We use the TR-BDF2 second-order scheme to discretize the Black-Scholes PDE \citep{lefloch2014tr}, using $N=1500$ time-steps, and vary the number of steps in the asset price dimension from $I=250$ to $I=4000$. The reference price is one obtained with $I=16000$, for the same $N$. It is close to the exact theoretical price, but it is different, since the number of time-steps is kept constant. The intent is to look at the convergence in the asset price dimension, not the overall convergence.

\subsection{Cubic vs. Sinh}
A uniform grid leads to largest error and oscillating convergence, because the accuracy depends strongly on the location of the critical point in the grid. The sinh stretching appear to be more accurate than the cubic stretching, convergence is somewhat more regular but still not of constant order for the same reasons as the uniform grid.
\begin{table}[h]
\centering{
	\caption{Absolute error in price $\times 10^5$ on a stretched grid. The reference price is obtained on a grid of $I=16000$ steps. \label{tbl:streching_noplacing}}
	\begin{tabular}{ccccccc}\toprule
$I$  & \multicolumn{3}{c}{$S=100$} & \multicolumn{3}{c}{$S=110$}\\\cmidrule(lr){2-4}\cmidrule(lr){5-7}
&  Uniform & Cubic & Sinh & Uniform & Cubic & Sinh\\\midrule
250 & 5002.3  & 1.3 & 256.9 & 5710.0 & 11.8 & 314.6 \\
500 & 74.0 & 387.8 & 71.7 & 89.1 & 434.5 & 73.7   \\
1000 & 1084.0 & 186.8 & 66.5 & 1223.9 & 209.9 & 76.5\\
2000 & 60.9 & 82.7 & 9.0 & 68.2 & 97.6 & 10.0\\\midrule
Reference Price & 2.31806 & 2.31735 & 2.31740 & 1.86342 &1.86263 & 1.86268 \\\bottomrule
	\end{tabular}
}
\end{table}

The choice $\alpha=1.5$ does not translate to exactly the same slope at the critical point for both transformations. The cubic transformation would require $\alpha=0.9$ to have the same slope. This partly explains the discrepancy in accuracy, with the reduced $\alpha$, the error with 500 points is significantly reduced to $137.8\times10^{-5}$.

\subsection{Placing vs. Deforming}
\subsubsection{Uniform}
\begin{figure}[h]
	\begin{center}
		%		\subfloat[][The 7-th sequence shows a small difference in favor of Kahan. Kahan is as accurate on the ordered sequences.]{\includegraphics[width=0.5\textwidth]{mean_kahan.pdf}}
		%		\subfloat[][Klein correction is much less accurate on the ordered sequences.]{\includegraphics[width=0.5\textwidth]{mean_klein_log.pdf}}
		\includegraphics[width=0.8\textwidth]{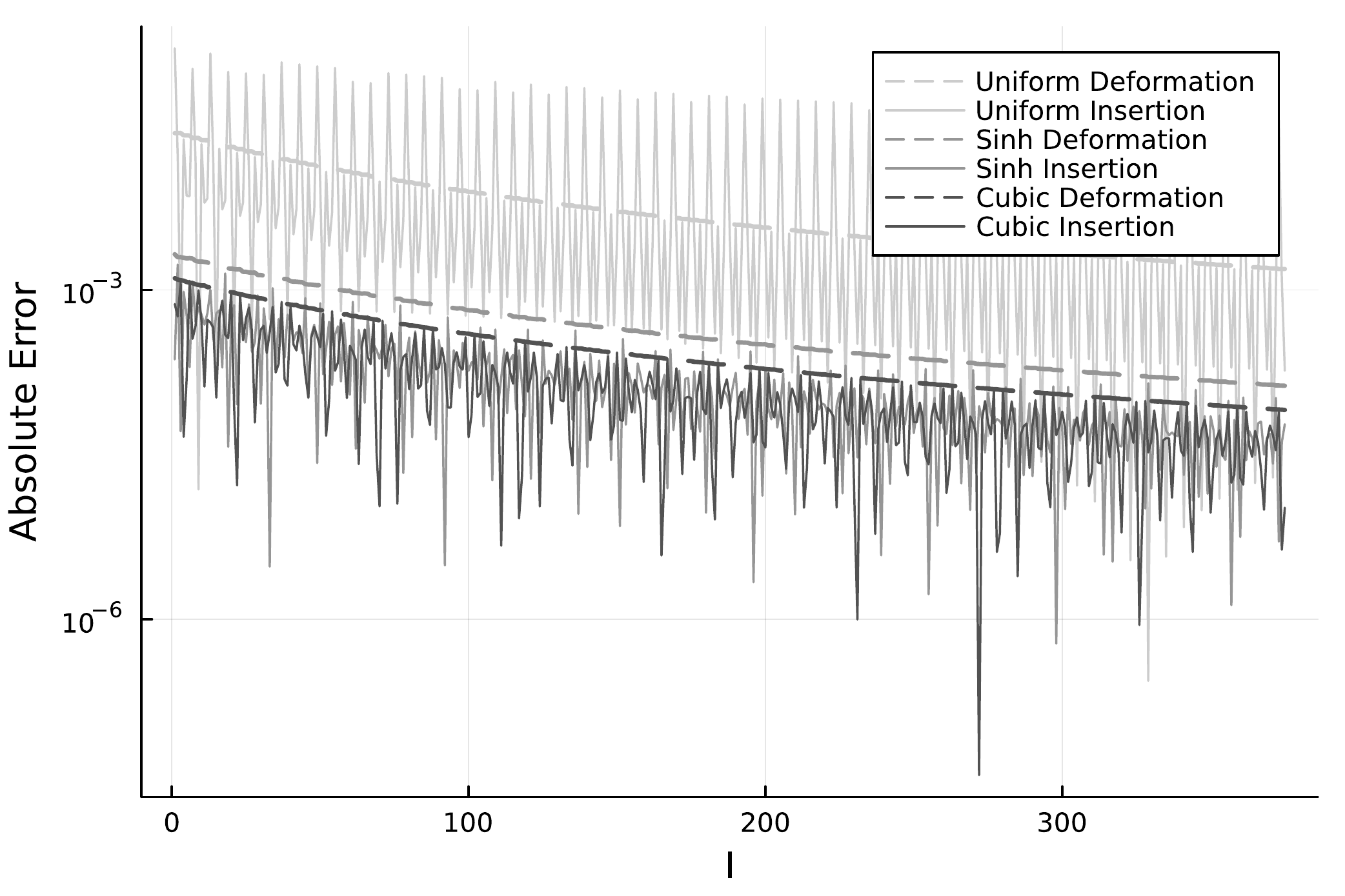}
	\end{center}
	\caption{Error in the price of a knock-out barrier option against the number of steps $I$ in the asset price dimension, for different kind of grids \label{fig:stretching-convergence}.}
\end{figure}

With the smooth grid deformation, the ratio of errors between doubling values $I$ is close to 4.0: the measured order of convergence order is close to two and stable (Table \ref{tbl:adjusted_uniform}). In contrast, the insertion of points does not lead to a smooth convergence.
\begin{table}[h]
	\centering{
		\caption{Absolute error in price $\times 10^5$ on an adjusted uniform grid. The reference price is obtained on a grid of $I=16000$ steps. \label{tbl:adjusted_uniform}}
		\begin{tabular}{ccccc}\toprule
			$I$  & \multicolumn{2}{c}{$S=100$} & \multicolumn{2}{c}{$S=110$}\\\cmidrule(lr){2-3}\cmidrule(lr){4-5}
			&  Deform & Insert & Deform & Insert \\ \midrule
			250 & 633.1  & 389.0 & 771.1 & 400.5 \\
			500 & 153.3 & 74.0 &  184.6 & 89.1 \\
			1000 & 38.2 & 98.5 & 46.5 & 108.8\\
			2000 & 9.4 & 60.9  & 11.3 & 68.1 \\\midrule
			Reference Price & 2.31736 & 2.31736 & 1.86264 &1.86263  \\\bottomrule
		\end{tabular}
	}
\end{table}
On this example, the insertion is less accurate than the deformation. This is slightly peculiar to the number of grid points and the location of the critical point. Figure \ref{fig:stretching-convergence} shows how much is the accuracy dependent on the grid details with the placing technique.
With the cubic or sinh stretching, the insertion is generally more accurate than the smooth deformation.

\subsubsection{Stretched}
Overall, the cubic stretching with insertion appear to be the most accurate on this problem (Table \ref{tbl:adjusted_stretched}). Figure \ref{fig:stretching-convergence} makes it however clear that the smooth deformation is preferable.
\begin{table}[h]
	\centering{
		\caption{Absolute error in price $\times 10^5$ of a knock-out barrier option, on an adjusted stretched grid. The reference price is obtained on a grid of $I=16000$ steps. \label{tbl:adjusted_stretched}}
		\begin{tabular}{ccccccccc}\toprule
			$I$  & \multicolumn{4}{c}{$S=100$} & \multicolumn{4}{c}{$S=110$}\\ \cmidrule(lr){2-5}\cmidrule(lr){6-9}
			& \multicolumn{2}{c}{Cubic} & \multicolumn{2}{c}{Sinh} & \multicolumn{2}{c}{Cubic} & \multicolumn{2}{c}{Sinh}\\\cmidrule(lr){2-3}\cmidrule(lr){4-5}  \cmidrule(lr){6-7} \cmidrule(lr){8-9}
			&  Deform  & Insert & Deform  &  Insert  &  Deform  & Insert & Deform &  Insert 			\\\midrule
			250 & 32.0  & 15.5 & 52.5 & 14.0 & 55.5 & 8.9 & 88.3 & 26.2  \\
			500 & 8.0 & 8.0 & 13.4 & 1.0 & 13.8 & 8.4 & 24.7 & 2.0 \\
			1000 & 2.0 & 1.9 & 3.3 & 2.3 & 3.8 & 2.4 & 5.5 & 4.1  \\
			2000 & 0.5 & 0.4 & 0.8 & 0.3 &  1.0 & 0.3 & 1.4 & 0.2\\\midrule
			Reference Price & 2.31736 & 2.31736 &2.31736 & 2.31736 & 1.86264 &1.86264 & 1.86264 &1.86264 \\\bottomrule
		\end{tabular}
	}
\end{table}
Placing the points result in a significant increase of the accuracy (along with a smooth convergence), this is particularly visible if we compare Table \ref{tbl:adjusted_stretched} with Table \ref{tbl:streching_noplacing}.

\subsection{Double Discrete Barrier Option}
In order to evaluate the difference between the stretching methods with multiple critical points, we price a double barrier put option with a strike price 102, maturity 0.5 year, up barrier level 110, down barrier level 90 under the Black-Scholes model with interest rate $r=10\%$, dividend yield $q=0\%$, volatility $\sigma=20\%$, underlying spot price $S=100$. We place the strike and the barriers in the middle of two grid points. The grid starts at $S_{\min}=54.57$ and ends at $S_{\max}=183.25$, which correspond to four standard deviations around the underlying spot price. We use $\alpha = 0.005(S_{\max}-S_{\min})$ to concentrate points around the three critical points.

The piecewise-cubic stretching is found to be as accurate as the Tavella Randall stretching (Table \ref{tbl:adjusted_stretched_double}). There is no obvious difference in terms of accuracy and convergence between $\mathcal{C}^2$ stretching and the cubic stretching.
\begin{table}[h]
	\centering{
		\caption{Absolute error in price $\times 10^5$ on an adjusted stretched grid of a double knock-out barrier with $\alpha=0.64$. The reference price of 0.600956 is obtained on a grid of $I=8000$ steps. \label{tbl:adjusted_stretched_double}}
		\begin{tabular}{ccccccccc}\toprule
			$I$  	& \multicolumn{2}{c}{Uniform} & \multicolumn{2}{c}{Piecewise Cubic} &  \multicolumn{2}{c}{Piecewise C2} & \multicolumn{2}{c}{Tavella Randall}\\\cmidrule(lr){2-3}\cmidrule(lr){4-5}  \cmidrule(lr){6-7}  \cmidrule(lr){8-9} 
			&  Deform  & Insert & Deform  &  Insert  & Deform  &  Insert &  Deform & Insert	\\\midrule
            50 & 6554.5 & 906.2 & 731.7 & 998.3 & 657.9 & 933.8 & 846.6 & 870.7  \\
			100 & 1658.7 & 542.9 & 207.2 &  233.4& 191.9 & 216.7 & 188.3 & 228.1  \\
			200 & 360.8 & 310.3& 49.8  &56.6  & 45.5 & 52.3 & 44.3 &47.5\\
			400 &  70.7& 44.5 &  10.8 & 11.7  & 9.7 & 10.7 & 8.1  &10.3  \\\bottomrule
		\end{tabular}
	}
\end{table}
Insertion is slightly worse than deformation on the stretched grid, the difference in accuracy is however relatively small. Insertion is more accurate on the uniform grid, but as evidenced in Figure \ref{fig:stretching-convergence}, this is highly dependent on the number of points used in the grid, i.e. where the inserted points fall in the initial uniform grid.

\subsection{American Option}
We consider an American put option contract of strike $K=100$ and maturity $T=1$, keeping otherwise the same Black-Scholes settings as in the previous numerical examples, and look at the convergence with number of steps $I$ in the asset price dimension for the different kinds of grid deformations.  In this problem, the second derivative of the solution is discontinuous around the exercise boundary and the first derivative is discontinuous at the strike price in the initial condition.
With a small $\alpha$ (relative to $S_{\max}-S_{\min}$), meaning a highly concentrated grid around the strike price, the error in the option price is almost the same as with a smoothly deformed uniform grid. The sinh stretching leads to a slightly higher error compared to the cubic stretching.
\begin{figure}[h]
	\begin{center}
				\subfloat[][$\alpha=1.50$.]{\includegraphics[width=0.5\textwidth]{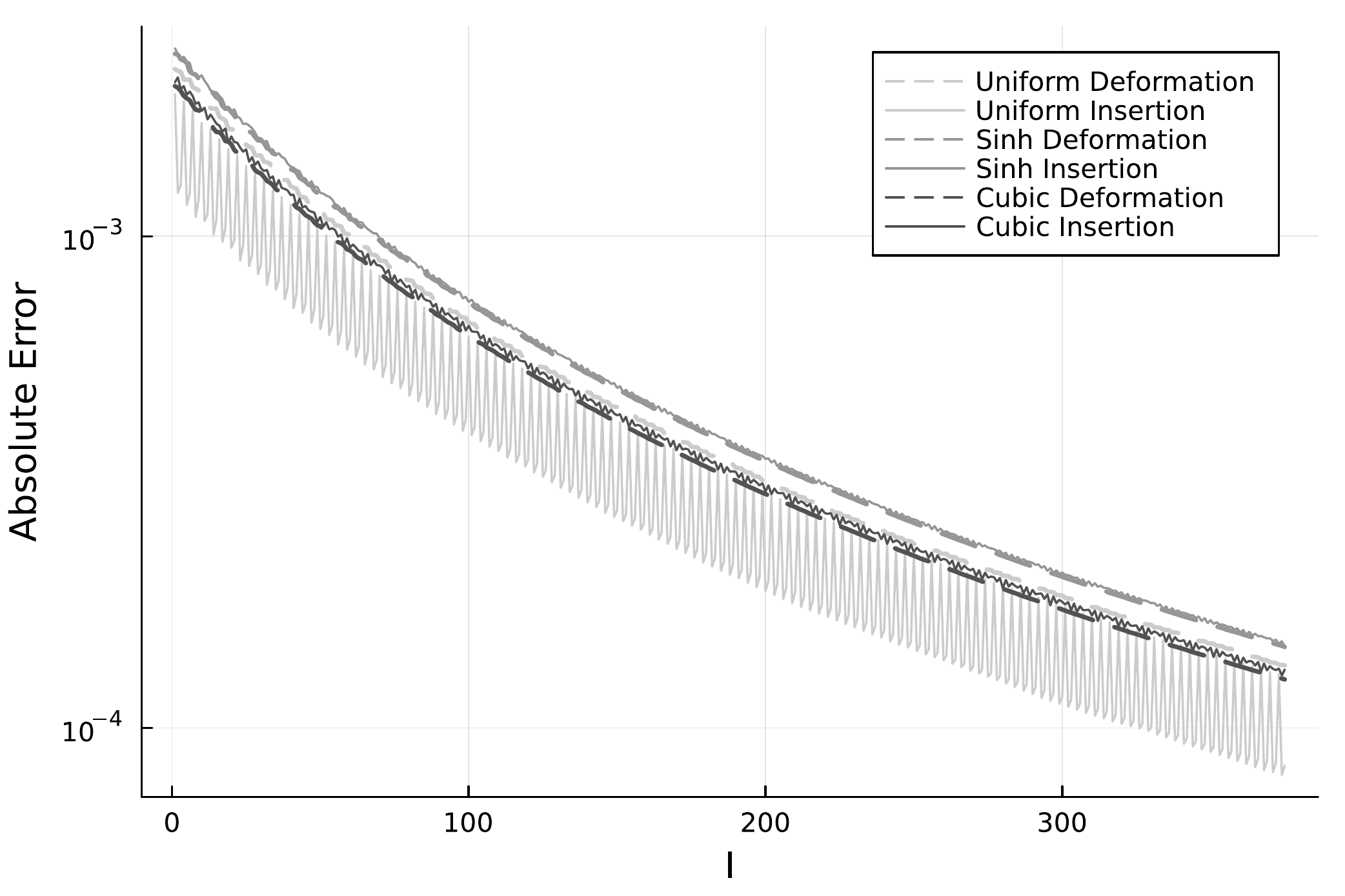}}
			\subfloat[][$\alpha=15.0$.]{\includegraphics[width=0.5\textwidth]{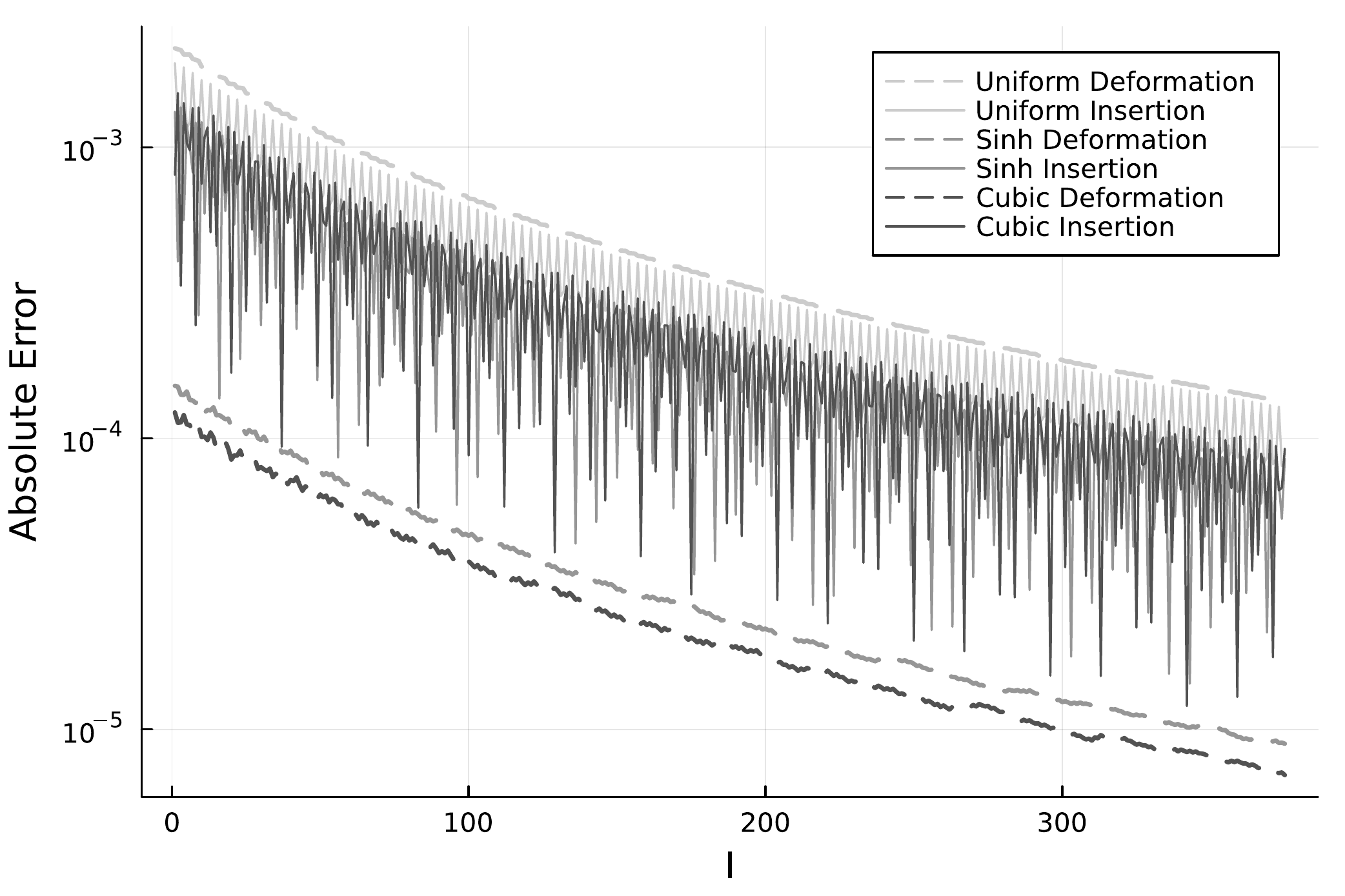}}
	\end{center}
	\caption{Error in the price of an American Put option, with different stretching around the point $K=100$. \label{fig:stretching_convergence_am}}
\end{figure}
With a larger $\alpha=15.0$, the transformations still concentrate points albeit less so than with the smaller $\alpha$, and the accuracy is much improved with the stretching. Insertion leads to clearly worse accuracy than a smooth deformation (Figure \ref{fig:stretching_convergence_am}).

\subsection{Continuous Barrier Option}
\subsubsection{Dirichlet boundary}
\begin{figure}[h!]
	\centering
	\begin{tikzpicture}
		\shade[shading=ball, ball color=black]  (0,7) circle (.2);
		\draw (0,8.2) node{$R$};
		\shade[shading=ball, ball color=white]  (2,7) circle (.2);
		\draw (2,7.8) node{$R$};
		\shade[shading=ball, ball color=black]  (4,7) circle (.2);
		\draw (4,8.2) node{$R$}; 
		\shade[shading=ball, ball color=white]  (6,7) circle (.2);
		\draw  (6,7.8) node{$R$};
		\shade[shading=ball, ball color=black]  (8,7) circle (.2);
		\draw  (8,8.2) node{$R$};
		\draw[ultra thick,-stealth] (0,8) -- (0,7.3);
		\draw[ultra thick,-stealth] (2,7.6) -- (2,7.3);
		\draw[ultra thick,-stealth] (4,8) -- (4,7.3);
		\draw[ultra thick,-stealth] (6,7.6) -- (6,7.3);
		\draw[ultra thick,-stealth] (8,8) -- (8,7.3);
		\draw[gray,,-stealth] (7.7,7) -- (6.3,7);
		\draw[gray,-stealth] (7.7,7) .. controls (5,7.7) .. (4.3,7);
		\draw[gray,-stealth] (6,7.3) .. controls (5.7,7.5) .. (4.3,7);
		\draw[gray,-stealth] (5.7,7) -- (4.3,7);
		\draw[style=dashed, line width=2pt] (0,5)--(8,5) (8.3,5) node[right]{$i_0$};
		\draw (8.3,7) node[right]{$i_0+1$};
		
		\shade[shading=ball, ball color=black]  (0,5) circle (.2);
		\draw (0,6.2) node{$R$};
		\shade[shading=ball, ball color=white]  (2,5) circle (.2);
		\draw (2,5.8) node{$R$};
		\shade[shading=ball, ball color=black]  (4,5) circle (.2);
		\draw (4,6.2) node{$R$}; 
		\shade[shading=ball, ball color=white]  (6,5) circle (.2);
		\draw  (6,5.8) node{$R$};
		\shade[shading=ball, ball color=black]  (8,5) circle (.2);
		\draw  (8,6.2) node{$R$};
		\draw[ultra thick,-stealth] (0,6) -- (0,5.3);
		\draw[ultra thick,-stealth] (2,5.6) -- (2,5.3);
		\draw[ultra thick,-stealth] (4,6) -- (4,5.3);
		\draw[ultra thick,-stealth] (6,5.6) -- (6,5.3);
		\draw[ultra thick,-stealth] (8,6) -- (8,5.3);
		\draw[gray,,-stealth] (7.7,5) -- (6.3,5);
		\draw[gray,-stealth] (7.7,5) .. controls (5,5.7) .. (4.3,5);
		\draw[gray,-stealth] (6,5.3) .. controls (5.7,5.5) .. (4.3,5);
		\draw[gray,-stealth] (5.7,5) -- (4.3,5);

		\draw[decoration={brace,mirror,raise=5pt},decorate,line width=2pt] (9.2,4.8) -- (9.2,7.2) (9.5,6) node[right]{\begin{tabular}{c} $A_{i,i-1}=0$,\\ $A_{i,i}=1$,\\ $A_{i,i+1}=0$. \end{tabular}};
		
		\shade[shading=ball, ball color=black]  (0,3) circle (.2);
		\shade[shading=ball, ball color=white]  (2,3) circle (.2);
		\shade[shading=ball, ball color=black]  (4,3) circle (.2);
		\shade[shading=ball, ball color=white]  (6,3) circle (.2);
		\shade[shading=ball, ball color=black]  (8,3) circle (.2);
		\draw  (8,4) node{$\max(S_{i_0-1}-K,0)$};
		\draw[ultra thick,-stealth] (8,3.8) -- (8,3.3);
		\draw[gray,,-stealth] (7.7,3) -- (6.3,3);
		\draw[gray,-stealth] (7.7,3) .. controls (5,3.7) .. (4.3,3);
		\draw[gray,-stealth] (6,3.3) .. controls (5.7,3.5) .. (4.3,3);
		\draw[gray,-stealth] (5.7,3) -- (4.3,3);
		\draw  (8.3,3) node[right]{$i_0-1$};

		\shade[shading=ball, ball color=black]  (0,1) circle (.2);
		\shade[shading=ball, ball color=white]  (2,1) circle (.2);
		\shade[shading=ball, ball color=black]  (4,1) circle (.2);
		\shade[shading=ball, ball color=white]  (6,1) circle (.2);
		\shade[shading=ball, ball color=black]  (8,1) circle (.2);
		\draw  (8,2) node{$\max(S_{i_0-2}-K,0)$};
		\draw[ultra thick,-stealth] (8,1.8) -- (8,1.3);
		\draw[gray,,-stealth] (7.7,1) -- (6.3,1);
		\draw[gray,-stealth] (7.7,1) .. controls (5,1.7) .. (4.3,1);
		\draw[gray,-stealth] (6,1.3) .. controls (5.7,1.5) .. (4.3,1);
		\draw[gray,-stealth] (5.7,1) -- (4.3,1);
		\draw  (8.3,1) node[right]{$i_0-2$};
	\end{tikzpicture}
	\caption{TR-BDF2 time-stepping of an up-and-out call option with rebate $R$ and strike $K$. The barrier falls on the grid at index $i_0$, of underlying asset value $S_{i_0}$. The matrix $A$ corresponds to the left hand-side of the linear tridiagonal system at a given time-step. Black balls are placed at actual grid times, white balls correspond to the intermediate stage.\label{fig:trbdf2_ko}}
\end{figure}
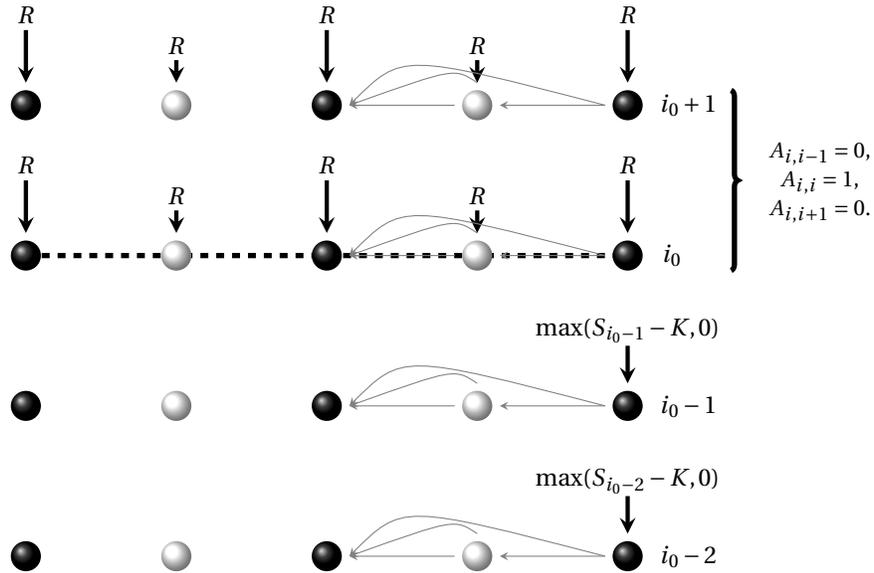

So far, our numerical examples were all around cases where the critical points are optimally located in the middle of two grid points. In the case of an option with a continuously monitored barrier, placing the barrier level exactly on the grid makes the boundary condition at the barrier simpler. For a knock-out, the value must be equal to the rebate at the barrier level (Figure \ref{fig:trbdf2_ko}), while for a knock-in we solve the PDE on two payoffs simultaneously; the vanilla option (result of the knock-in) and the knock-in option itself where the boundary condition states that at and above the barrier the value must be the one of the European option.

\subsubsection{Ghost point}
A more general technique, which works when the barrier level is not on the grid, consists in the use of ghost point (on the grid) above the barrier, which ensures that the value at the barrier is exactly zero. The trapezoidal stage of TR-BDF2 consists in an explicit Euler step followed by an implicit Euler step. Let $i_0$ be the index such that $S_{i_0-1} < B < S_{i_0}$, and let us consider the timestep from the time $t_j$ to the time $t_{j-1}$. The value at the ghost point must be such that $V(B,t_{j-1}) = R$. \citet[p. 1209--1210]{wilmott2013paul} suggests that a linear interpolation is good enough to determine the value  $G$ at the ghost point. This leads to
\begin{equation}
	 \frac{B- S_{i_0-1}}{S_{i_0}-S_{i_0-1}}G +  \frac{S_{i_0}-B}{S_{i_0}-S_{i_0-1}} V(t_j, S_{i_0-1}) = R\,, \label{eqn:ghost_point_linear_explicit}
\end{equation}
or equivalently
\begin{equation*}
  V(t_j, S_{i_0}) =  G = \frac{(S_{i_0}-S_{i_0-1})R - (S_{i_0}-B)V(t_j, S_{i_0-1})}{B- S_{i_0-1}}\,. 
\end{equation*}
Alternatively, a three-points Lagrange interpolation would lead to
\begin{align}
 G  \frac{(B-S_{i_0-2})(B-S_{i_0-1})}{(S_{i_0}-S_{i_0-2})(S_{i_0}-S_{i_0-1})} +  V(t_j, S_{i_0-1}) \frac{(B-S_{i_0-2})(B-S_{i_0})}{(S_{i_0-1}-S_{i_0-2})(S_{i_0-1}-S_{i_0})} \nonumber\\ +  V(t_j, S_{i_0-2}) \frac{(B-S_{i_0})(B-S_{i_0-1})}{(S_{i_0-2}-S_{i_0})(S_{i_0-2}-S_{i_0-1})} = R\,, \label{eqn:ghost_point_lagrange_explicit}
\end{align}
from which we may deduce $G$. Before calculating the right hand side of the linear system, we thus override the value of $V(t_j, S_{i_0})$ with the one obtained through Equation \ref{eqn:ghost_point_linear_explicit} or Equation \ref{eqn:ghost_point_lagrange_explicit} and use a Dirichlet boundary condition at $i_0$.

Similar care needs to be taken for the implicit part of the trapezoidal stage. In this case, we enforce Equation \ref{eqn:ghost_point_linear_explicit} at $t_{j-1}$ instead of $t_j$. The left hand side matrix $A$ of the linear system is thus modified such that 
\begin{align*}
	A_{i_0,i_0} = \frac{B- S_{i_0-1}}{S_{i_0}-S_{i_0-1}}\,,&\quad A_{i_0,i_0-1} = \frac{S_{i_0}-B}{S_{i_0}-S_{i_0-1}}\,,\quad A_{i_0,i_0+1}=0 \,.
\end{align*}
The three-points interpolation leads to an additional term $A_{i_0,i_0-2}$, which may be removed using a linear combination of the system at row $i_0-1$ in order to keep a tridiagonal system.

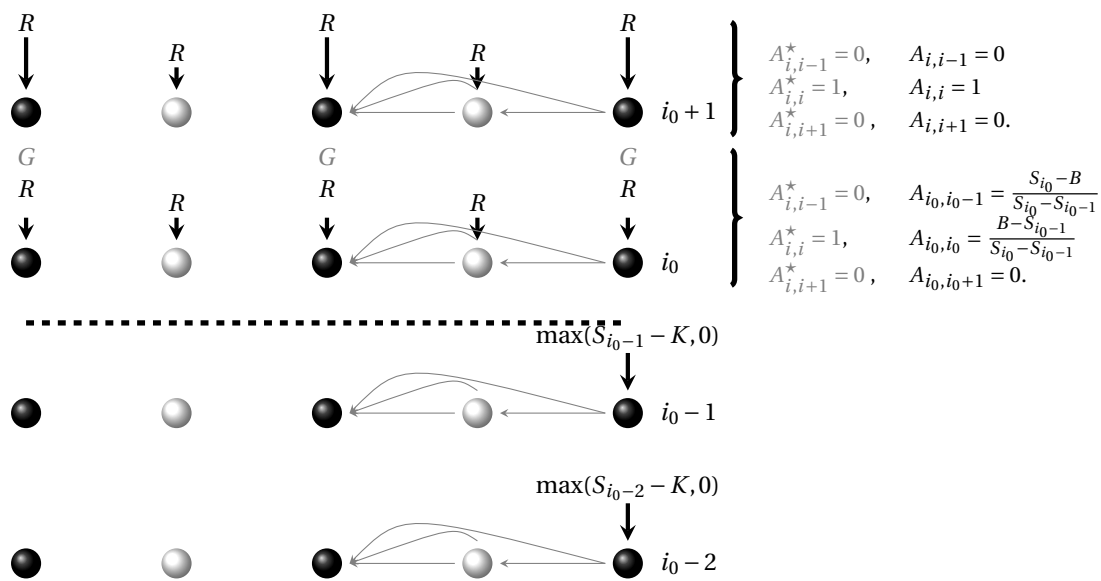
\begin{figure}[h!]
	\centering
	\begin{tikzpicture}
		\shade[shading=ball, ball color=black]  (0,7) circle (.2);
		\draw (0,8.2) node{$R$};
		\shade[shading=ball, ball color=white]  (2,7) circle (.2);
		\draw (2,7.8) node{$R$};
		\shade[shading=ball, ball color=black]  (4,7) circle (.2);
		\draw (4,8.2) node{$R$}; 
		\shade[shading=ball, ball color=white]  (6,7) circle (.2);
		\draw  (6,7.8) node{$R$};
		\shade[shading=ball, ball color=black]  (8,7) circle (.2);
		\draw  (8,8.2) node{$R$};
		\draw[ultra thick,-stealth] (0,8) -- (0,7.3);
		\draw[ultra thick,-stealth] (2,7.6) -- (2,7.3);
		\draw[ultra thick,-stealth] (4,8) -- (4,7.3);
		\draw[ultra thick,-stealth] (6,7.6) -- (6,7.3);
		\draw[ultra thick,-stealth] (8,8) -- (8,7.3);
		\draw[gray,,-stealth] (7.7,7) -- (6.3,7);
		\draw[gray,-stealth] (7.7,7) .. controls (5,7.7) .. (4.3,7);
		\draw[gray,-stealth] (6,7.3) .. controls (5.7,7.5) .. (4.3,7);
		\draw[gray,-stealth] (5.7,7) -- (4.3,7);
		\draw[style=dashed, line width=2pt] (0,4.2)--(8,4.2) ;
		\draw (8.3,5) node[right]{$i_0$};
		\draw (8.3,7) node[right]{$i_0+1$};
		
		\shade[shading=ball, ball color=black]  (0,5) circle (.2);
		\draw (0,6.2)  node[text width=0.5cm, align=center]{$\textcolor{gray}{G}$ \linebreak $R$};
		\shade[shading=ball, ball color=white]  (2,5) circle (.2);
		\draw (2,5.8) node{$R$};
		\shade[shading=ball, ball color=black]  (4,5) circle (.2);
		\draw (4,6.2)  node[text width=0.5cm, align=center]{$\textcolor{gray}{G}$ \linebreak $R$}; 
		\shade[shading=ball, ball color=white]  (6,5) circle (.2);
		\draw  (6,5.8) node{$R$};
		\shade[shading=ball, ball color=black]  (8,5) circle (.2);
		\draw  (8,6.2) node[text width=0.5cm, align=center]{$\textcolor{gray}{G}$ \linebreak $R$};
		\draw[ultra thick,-stealth] (0,5.6) -- (0,5.3);
		\draw[ultra thick,-stealth] (2,5.6) -- (2,5.3);
		\draw[ultra thick,-stealth] (4,5.6) -- (4,5.3);
		\draw[ultra thick,-stealth] (6,5.6) -- (6,5.3);
		\draw[ultra thick,-stealth] (8,5.6) -- (8,5.3);
		\draw[gray,,-stealth] (7.7,5) -- (6.3,5);
		\draw[gray,-stealth] (7.7,5) .. controls (5,5.7) .. (4.3,5);
		\draw[gray,-stealth] (6,5.3) .. controls (5.7,5.5) .. (4.3,5);
		\draw[gray,-stealth] (5.7,5) -- (4.3,5);

		\draw[decoration={brace,mirror,raise=5pt},decorate,line width=2pt] (9.2,6.7) -- (9.2,8.2) (9.5,7.3) node[right]{\begin{tabular}{ll} 	
				$\textcolor{gray}{A^\star_{i,i-1}=0}$, &  $A_{i,i-1} = 0$\\ $\textcolor{gray}{A^\star_{i,i}=1}$,& 	$A_{i,i} = 1$\\ $\textcolor{gray}{A^\star_{i,i+1}=0}$ ,&  $A_{i,i+1}=0$. \end{tabular}};
		
		\draw[decoration={brace,mirror,raise=5pt},decorate,line width=2pt] (9.2,4.7) -- (9.2,6.5) (9.5,5.4) node[right]{\begin{tabular}{ll} 
				$\textcolor{gray}{A^\star_{i,i-1}=0}$, &  $A_{i_0,i_0-1} = \frac{S_{i_0}-B}{S_{i_0}-S_{i_0-1}}$\\ $\textcolor{gray}{A^\star_{i,i}=1}$,& 	$A_{i_0,i_0} = \frac{B- S_{i_0-1}}{S_{i_0}-S_{i_0-1}}$\\ $\textcolor{gray}{A^\star_{i,i+1}=0}$ ,&  $A_{i_0,i_0+1}=0$. \end{tabular}};
		
		\shade[shading=ball, ball color=black]  (0,3) circle (.2);
		\shade[shading=ball, ball color=white]  (2,3) circle (.2);
		\shade[shading=ball, ball color=black]  (4,3) circle (.2);
		\shade[shading=ball, ball color=white]  (6,3) circle (.2);
		\shade[shading=ball, ball color=black]  (8,3) circle (.2);
		\draw  (8,4) node{$\max(S_{i_0-1}-K,0)$};
		\draw[ultra thick,-stealth] (8,3.8) -- (8,3.3);
		\draw[gray,,-stealth] (7.7,3) -- (6.3,3);
		\draw[gray,-stealth] (7.7,3) .. controls (5,3.7) .. (4.3,3);
		\draw[gray,-stealth] (6,3.3) .. controls (5.7,3.5) .. (4.3,3);
		\draw[gray,-stealth] (5.7,3) -- (4.3,3);
		\draw  (8.3,3) node[right]{$i_0-1$};

		\shade[shading=ball, ball color=black]  (0,1) circle (.2);
		\shade[shading=ball, ball color=white]  (2,1) circle (.2);
		\shade[shading=ball, ball color=black]  (4,1) circle (.2);
		\shade[shading=ball, ball color=white]  (6,1) circle (.2);
		\shade[shading=ball, ball color=black]  (8,1) circle (.2);
		\draw  (8,2) node{$\max(S_{i_0-2}-K,0)$};
		\draw[ultra thick,-stealth] (8,1.8) -- (8,1.3);
		\draw[gray,,-stealth] (7.7,1) -- (6.3,1);
		\draw[gray,-stealth] (7.7,1) .. controls (5,1.7) .. (4.3,1);
		\draw[gray,-stealth] (6,1.3) .. controls (5.7,1.5) .. (4.3,1);
		\draw[gray,-stealth] (5.7,1) -- (4.3,1);
		\draw  (8.3,1) node[right]{$i_0-2$};
	\end{tikzpicture}
	\caption{TR-BDF2 time-stepping of an up-and-out call option with rebate $R$ and strike $K$. The barrier falls in between grid points. The matrix $A$ corresponds to the left hand-side of the linear tridiagonal system at a given time-step. $A^\star$ and gray values correspond to the right hand-side of the trapezoidal step. Black balls are placed at actual grid times, white balls correspond to the intermediate stage.\label{fig:trbdf2_ko_ghost}}
\end{figure}

The BDF2 stage may reuse the same left hand side matrix. The overall technique is summarized in Figure \ref{fig:trbdf2_ko_ghost}.

We evaluate the accuracy of the three-points Lagrange interpolation to determine the boundary condition in the context of a double knock out option of strike $K=100$, maturity $T=1$ year, lower barrier $L=90$, upper barrier $U=160$ on an asset of spot price $S=95$, with interest rate $r=10\%$ and Black-Scholes volatility $\sigma=25\%$ in Table \ref{tbl:continuous_double_ghost}. The grid is uniform composed of the same number of time-steps and space-steps, truncated at or just above the barrier levels.
\begin{table}[h!]
	\centering{
		\caption{Absolute error in price $\times 10^5$ for a continuously monitored double knock-out barrier. The reference price of 3.460714 is obtained on a grid of $I=4000$ space-steps and the same number of time-steps. \label{tbl:continuous_double_ghost}}
		\begin{tabular}{cccc}\toprule
			$I$ & Ghost (Linear) & Ghost (3-points) & On Grid 		\\\midrule
			20 &  3801.7 & 825.9 &  652.2  \\
			40 & 899.1 &  178.9 & 125.7\\
			80 & 230.0 & 28.9 & 27.1 \\
			160 &  52.3 & 11.1 & 9.6  \\\bottomrule
		\end{tabular}
	}
\end{table}
On this problem, the linear interpolation results in an important loss of accuracy, while the three points interpolation leads to nearly the same accuracy as the case where the barrier is on the grid. Finally, if we move slightly the grid, the accuracy of the ghost point technique does not vary. Those conclusions stand for different values of the volatility or  of the interest rate.

The ghost point technique is more involved to implement, and makes more sense in the context of a time-varying barrier, such as, for example an exponential barrier. On such contracts, the three-points approximation is not necessarily more accurate, as suggested in \citep{wilmott2013paul}. In practice however, unless the underlying variable has been transformed in time, the barrier is flat (or piecewise-flat) in financial derivative contracts, because of operational constraints.

% Space steps   Time steps   Value  Error (Lin)   Value Error (Quad)   Value Error (On Grid)

\subsubsection{Stretched}
In Table \ref{tbl:continuous_double_stretching} we look at the same continuous double knock-out option as in the previous section, this time, using a deformation to place the strike in the middle and the barriers on the grid and eventually stretch the grid (using $\alpha=0.01(U-L)=0.7$) to be more dense near those three critical points.
\begin{table}[h]
	\centering{
		\caption{Absolute error in price $\times 10^5$ for a continuously monitored double knock-out barrier on a deformed grid such that the barriers fall on the grid exactly. The reference price of 3.460714 is obtained on a grid of $I=4000$ space-steps and the same number of time-steps. \label{tbl:continuous_double_stretching}}
		\begin{tabular}{cccc}\toprule
			$I$ & Uniform & Piecewise Cubic & Tavella Randall 		\\\midrule
			20 &  772.5 & 6204.4 &18065.3  \\
			40 & 173.2 &  1440.7 & 4260.2\\
			80 & 38.0 & 326.5 &1010.4 \\
			160 &  9.6 & 80.2 &238.6   \\\bottomrule
		\end{tabular}
	}
\end{table}
We find that concentrating points close to the barriers does not increase the accuracy. On the contrary, we notice a significant drop in accuracy with the Tavella-Randall or the cubic stretchings on this problem. This also holds if we increase the density only around the strike price as well as if we use a ghost point instead of placing the barriers.

%\subsubsection{Insertion}
%Insertion vs displacement?
%In order to place the barrier level on the grid, we evaluate two different ways:
%\begin{itemize}
%	\item Insertion: we simply insert the barrier level. We may wonder if having two points too close by could be a problem. It certainly would be for an explicit finite difference scheme.
%	\item Displacement: we move slightly the closest point to the barrier at the barrier level.
%\end{itemize}

%Error as a function of the number of steps in the asset price dimension. Eventually barrier will fall close to a grid point. => no significant difference in accuracy between the two choices for the TR-BDF2 scheme.

%TODO : evaluate inverse lin (book) interpolation vs tavella approach for placement of Sinh. Could be Java ExactStretching impl.
%TODO: evaluate lagrange 3 points vs cublic spline for gamma in book. Lagrange 3 points seems good enough. implies gamma is contant and jumps? 

\section{Conclusion}
Inserting points such that the critical points fall in the middle of two grid points increases the accuracy compared to a raw uniform grid in most situations. It is also effective on stretched grids. A smooth deformation via a cubic spline is however almost always preferable, and leads to a smooth convergence. Furthermore, the smooth deformation always enhances significantly the accuracy when applied on top of a preexisting grid-stretching. 

In terms of stretching, the simple cubic transformation is found to be at least as accurate as the hyperbolic sine transformation, while using less computational resources. This is even more relevant when the problem involves multiple critical points.

Finally, stretching is very effective on discrete barrier or American options but may be detrimental on occasion, for example when applied to a continuously monitored double barrier option.

\funding{This research received no external funding.}
\conflictsofinterest{The authors declare no conflict of interest.}
%\supplementary{Source code illustrating the examples in this paper is available at \url{https://github.com/jherekhealy/MonteCarloMeanVarianceExamples.jl}}
\externalbibliography{yes}
\bibliography{stretching.bib}
\appendixtitles{no}
\end{document}